\documentstyle[12pt]{article}
\newfont{\frak}{eufm10 scaled\magstep1}
\newfont{\Bbb}{msbm10 scaled\magstep1}                                     
\let\al\alpha
\let\be\beta
\let\ga\gamma
\let\de\delta

\let\ka\kappa
\let\la\lambda
\let\Ga\Gamma
\let\De\Delta
\newcommand{\ZZ}{\hbox{\Bbb Z}}
\newcommand{\CC}{\hbox{\Bbb C}}
\newcommand{\RR}{\hbox{\Bbb R}}
\newcommand{\gog}{\hbox{\frak g}}
\newcommand{\FSA}{{\cal A}}

\newcommand{\FSD}{{\cal D}}
\newcommand{\FSG}{{\cal G}}
\newcommand{\FSO}{{\cal O}}
\newcommand{\FSU}{{\cal U}}
\let\pa\partial
\let\iy\infty
\newcommand{\halmos}{\hbox{\vrule height0.28cm width0.01cm\vbox{\hrule height
 0.01cm width0.3cm \vskip0.26cm \hrule height 0.01cm width0.3cm}\vrule
 height0.28cm width 0.01cm}}
\newcommand{\eof}{\quad\halmos}
\newcommand{\tfrac}[2]{{\textstyle\frac{#1}{#2}}}
\newcommand{\const}{{\rm const.\,}}
\newcommand{\Proof}{\noindent{\bf Proof}\ \ }

\newtheorem{theorem}{Theorem}[section]
\newtheorem{prop}[theorem]{Proposition}

\newtheorem{lemma}[theorem]{Lemma}
\newtheorem{Remark}[theorem]{Remark}
\newenvironment{remark}{\begin{Remark}\rm\ }{\end{Remark}}

\begin{document}
\vskip 50mm
\title{Gauss hypergeometric function and \\
quadratic $R$-matrix algebras }

\vskip 45mm
\author{Tom H. Koornwinder and Vadim B. Kuznetsov${}\sp {1}$}
\footnotetext[1]{This author was supported by the
Netherlands Organisation for Scientific Research
(NWO) under the Project \# 611--306--540.}

\vskip 10mm
\maketitle

{\it Department of Mathematics and Computer Science, 
University of Amsterdam,
\par Plantage Muidergracht 24, 
1018 TV Amsterdam, The Netherlands
\par e-mail: {\tt thk@fwi.uva.nl}, {\tt vadim@fwi.uva.nl}}

\vskip 10mm
This paper is dedicated to L. D. Faddeev on the occasion
of his sixtieth birthday. 

\vskip 20mm
\begin{abstract}\noindent
We consider  representations
of quadratic $R$-matrix algebras
by means of certain first order ordinary differential operators.
These operators turn out to act as parameter shifting operators on
the Gauss hypergeometric function and its limit cases and on classical 
orthogonal polynomials. The relationship with W. Miller's treatment of
Lie algebras of first order differential operators will be discussed.

\end{abstract}

\par classical orthogonal polynomials, 
recurrence relations
\vskip 10mm
{\bf AMS classification:} 33C05, 33C35, 58F07 

\vskip 10mm
\pagebreak
\section{Introduction}
\setcounter{equation}{0}

The modern approach to finite-dimensional integrable systems 
uses the language of the 
representations of $R$-matrix algebras 
\cite{ft87,ks82,ku90,rs87,sk84,sk91}. 
There are two quadratic $R$-matrix algebras appearing in
the quantum inverse scattering method (QISM). We will call them
QISM I and QISM II, respectively. We restrict ourselves to the 
simplest case of a 2-dimensional auxiliary space and a
rational $4 \times 4$ $R$-matrix of the form
\begin{equation}
R(u)=u+\ka P, \quad P=\left(\matrix{
1&0&0&0\cr 0&0&1&0\cr 0&1&0&0\cr 0&0&0&1}\right),\qquad \ka,u\in\CC.
\label{1}
\end{equation} 
Consider a $2\times 2$ matrix
\begin{equation}
\left(\matrix{A(u)&B(u)\cr C(u)&D(u)}\right)
\label{2}\end{equation}
with a priori non-commuting entries depending on a so-called
{\em spectral parameter} $u$ which is arbitrary complex.
The matrix (\ref{2}) is denoted $T(u)$ in the QISM I case and $U(u)$ in the
QISM II case.
The {\em QISM I algebra} or {\em $T$-algebra}
is then defined as the algebra generated by all matrix elements
of $T(u)$ for all complex values of $u$ subject to the following quadratic
relation on $T(u)$ (cf.\ \cite{ba82,ks82,tf79}).
\begin{equation}
R(u-v)T^{(1)}(u)T^{(2)}(v)=T^{(2)}(v)T^{(1)}(u)R(u-v),\quad u,v\in\CC.
\label{3}\end{equation} 
Here we use the notation
$T^{(1)}(u)=T(u)\otimes I$, \,$T^{(2)}(v)=I\otimes T(v)$.
The {\em QISM II algebra} or {\em $U$-algebra} is the algebra generated by the
matrix elements of $U(u)$ for all $u$ subject to a quadratic
relation involving two $R$-matrices \cite{sk88}:
\begin{equation}
R(u-v)U^{(1)}(u)R(u+v-1)U^{(2)}(v)=
U^{(2)}(v)R(u+v-1)U^{(1)}(u)R(u-v),\quad u,v\in\CC.
\label{4}\end{equation}
{}From now on we assume for both types of algebras that $\ka=1$. For $\ka\ne0$
in (\ref{1}) this means no loss of generality.

In the present article we construct representations of very simple type
($L$-operators of rank 1) of both the $T$- and the $U$-algebra.
In the QISM II case we require moreover a certain symmetry property (unitarity)
for the $L$-operator.
We will consider $L$-operators (\ref{2}) for which certain matrix elements
will be realized as first order ordinary differential operators acting
as parameter shifting operators on the Gauss hypergeometric function
and its limit cases. Specialization then yields shift operator actions on
classical orthogonal polynomials.
For the QISM I case and for some of the QISM II cases
we will point out a close connection of our results
with the Infeld-Hull \cite{infeld} factorization method for second order
differential equations and with Miller's \cite{mi68} treatment of
Lie algebras of first order differential operators acting as shift operators
on special functions. For the most general QISM II cases, we consider,
the connection
with Lie algebras of first order differential operators is no longer valid.
But then, instead, there is a connection with an action by
differential operators (cf.\ Miller \cite{mi68})
of the universal enveloping algebra
of the Lie algebra ${\rm e}(3)$.

Our operators will act on special functions $F(u)$
which appear for each $u\in\CC$ as a solution of the equation
\begin{equation}
C(u)\,F(u)=0,
\label{5}\end{equation} 
i.e., as functions annihilated by one of the two off-diagonal
elements (always chosen to be $C(u)$)
of an $L$-operator. The operators $A(u)$ and $D(u)$,
for the QISM I algebra, and $A(u)$ and $-A(-u)$, for the QISM II
algebra, then give the shifting of the parameter $u$ by $\pm 1$,
respectively:
\begin{equation}
A(u)F(u)=\Delta_-(u-\tfrac12)F(u-1),\quad 
D(u)F(u)=\Delta_+(u+\tfrac12)F(u+1),
\label{6}\end{equation}
for the $T$-algebra, and
\begin{equation}
A(u)F(u)=\Delta_-(u-\tfrac12)F(u-1),\quad 
-A(-u)F(u)=\Delta_+(u+\tfrac12)F(u+1),
\label{7}\end{equation}
for the $U$-algebra. Here the $\Delta_\pm(u)$ are certain scalars
depending on $u$ which factorize the quantum determinant $\Delta(u)$ 
of an $L$-operator:
\begin{equation}
\Delta(u)=\Delta_+(u)\Delta_-(u).
\label{8}\end{equation}
The quantum determinant of a $T$- or $U$-algebra is a certain quadratic
expression in the
generators with the property that it is the generating function for the
center of the algebra.
So, in an irreducible representation it is,
under suitable assumptions, scalar for each $u$.

The structure of the paper is as follows. In Section 2 we 
give further 
properties of both types of algebras. Section 3 contains a representative
collection of differential recurrence relations for
special functions for which we can give an interpretation in terms
of $L$-operators satisfying (\ref{3}) or (\ref{4}).
Section 4 clarifies the connection of our approach with the
factorization method and with Miller's Lie algebra approach.   
Section 5 deals with the  
simplest (rank 1) $L$-operators for the QISM I algebra and 
with the corresponding shifting formulas for Gauss
hypergeometric functions, etc. In Section 6 we study
rank 1 $L$-operators for the QISM II algebra and the corresponding
differential recurrence relations.
In the final Section 7 we make some concluding 
remarks on possible applications.

Throughout we use notation like $\pa_x$ for the ordinary derivative $d/dx$ or
the partial derivative $\pa/\pa x$ with respect to $x$.

\section{More about quadratic $R$-matrix algebras}
\setcounter{equation}{0}
In the QISM II algebra case we will always add the following relations
(symmetry property when changing the sign of $u$):
\begin{eqnarray}
-A(-u)&=&D(u)-(A(u)+D(u))/(2u+1),\nonumber\\
-D(-u)&=&A(u)-(A(u)+D(u))/(2u+1), 
\nonumber\\
B(-u)&=&B(u),\quad C(-u)=C(u).
\label{141}\end{eqnarray}
Note that the second equality is implied by the first.
The equations (\ref{141}) can be rephrased as the unitarity property
$U^{-1}(-u)\sim U(u)$ (\cite{sk88}).

The {\em quantum determinant} of a $T$- or $U$-algebra is defined as follows.
\begin{eqnarray}
\Delta(u)
&=&A(u-\tfrac12)D(u+\tfrac12)-C(u-\tfrac12)B(u+\tfrac12)
\nonumber\\
&=&D(u-\tfrac12)A(u+\tfrac12)-B(u-\tfrac12)C(u+\tfrac12)
\nonumber\\
&=&D(u+\tfrac12)A(u-\tfrac12)-C(u+\tfrac12)B(u-\tfrac12)
\nonumber\\
&=&A(u+\tfrac12)D(u-\tfrac12)-B(u+\tfrac12)C(u-\tfrac12),
\label{9}\end{eqnarray}  
for the $T$-algebra, and
\begin{eqnarray}
\Delta(u)
&=&-D(-u+\tfrac12)D(u+\tfrac12)-C(u-\tfrac12)B(u+\tfrac12)
\nonumber\\
&=&-A(-u+\tfrac12)A(u+\tfrac12)-B(u-\tfrac12)C(u+\tfrac12)
\nonumber\\
&=&-D(u+\tfrac12)D(-u+\tfrac12)-C(u+\tfrac12)B(u-\tfrac12)
\nonumber\\
&=&-A(u+\tfrac12)A(-u+\tfrac12)-B(u+\tfrac12)C(u-\tfrac12),
\label{10}\end{eqnarray}  
for the $U$-algebra. The quantum determinant is the generating
function for the center of both types of algebras \cite{ks82,sk88}. 

Relation (\ref{3}) resp.\ (\ref{4}) can be rewritten in the following
extended form as commutators between the algebra generators
$A(u),B(u),C(u)$, and $D(u)$.
\begin{eqnarray}
{[A,A]}&=&[B,B]=[C,C]=[D,D]=0,\label{21}\\
{[A,B]}&=&- (AB-{\tilde {AB}})/(u-v),\label{22}\\
{[B,A]}&=&- (BA-{\tilde {BA}})/(u-v),\label{23}\\
{[A,C]}&=&- (CA-{\tilde {CA}})/(u-v),\label{24}\\
{[C,A]}&=&- (AC-{\tilde {AC}})/(u-v),\label{25}\\
{[B,D]}&=&- (DB-{\tilde {DB}})/(u-v),\label{26}\\
{[D,B]}&=&- (BD-{\tilde {BD}})/(u-v),\label{27}\\
{[D,C]}&=&- (DC-{\tilde {DC}})/(u-v),\label{28}\\
{[C,D]}&=&- (CD-{\tilde {CD}})/(u-v),\label{29}\\
{[A,D]}&=&- (CB-{\tilde {CB}})/(u-v),\label{210}\\
{[D,A]}&=&- (BC-{\tilde {BC}})/(u-v),\label{211}\\
{[B,C]}&=&- (DA-{\tilde {DA}})/(u-v),\label{212}\\
{[C,B]}&=&- (AD-{\tilde {AD}})/(u-v),\label{213}
\end{eqnarray}
for the $T$-algebra, and
\begin{eqnarray}
{[B,B]}&=&[C,C]=0,\label{214}\\
{[A,A]}&=&- (BC-{\tilde {BC}})/(u+v),\label{215}\\
{[D,D]}&=&- (CB-{\tilde {CB}})/(u+v),\label{216}\\
{[A,B]}&=&- (AB-{\tilde {AB}})/(u-v)-
(AB+BD)/(u+v-1)\nonumber\\
&&-(AB+BD-{\tilde {AB}}-{\tilde {BD}})/
(u-v)/(u+v-1),\label{217}\\
{[B,A]}&=&- (BA-{\tilde {BA}})/(u-v)+
({\tilde {AB}}+{\tilde {BD}})/(u+v-1),\label{218}\\
{[A,C]}&=&- (CA-{\tilde {CA}})/(u-v)+
({\tilde {CA}}+{\tilde {DC}})/(u+v-1)\nonumber\\
&&-(CA+DC-{\tilde {CA}}-{\tilde {DC}})/
(u-v)/(u+v-1),\label{219}\\
{[C,A]}&=&- (AC-{\tilde {AC}})/(u-v)-
(CA+DC)/(u+v-1),\label{220}\\
{[D,B]}&=&- (BD-{\tilde {BD}})/(u-v)+
({\tilde {AB}}+{\tilde {BD}})/(u+v-1)\nonumber\\
&&-(AB+BD-{\tilde {AB}}-{\tilde {BD}})/
(u-v)/(u+v-1),\label{221}\\
{[B,D]}&=&- (DB-{\tilde {DB}})/(u-v)-
(AB+BD)/(u+v-1),\label{222}\\
{[D,C]}&=&- (DC-{\tilde {DC}})/(u-v)-
(CA+DC)/(u+v-1)\nonumber\\
&&-(CA+DC-{\tilde {CA}}-{\tilde {DC}})/
(u-v)/(u+v-1),\label{223}\\
{[C,D]}&=&- (CD-{\tilde {CD}})/(u-v)+
({\tilde {CA}}+{\tilde {DC}})/(u+v-1),\label{224}\\
{[A,D]}&=&- (CB-{\tilde {CB}})(u+v+1)/(u^2-v^2),\label{225}\\
{[D,A]}&=&- (BC-{\tilde {BC}})(u+v+1)/(u^2-v^2),\label{226}\\
{[B,C]}&=&- (DA-{\tilde {DA}})(u+v-1)/(u^2-v^2)\nonumber\\
&&- (AA-{\tilde {DD}})/(u+v),\label{227}\\
{[C,B]}&=&- (AD-{\tilde {AD}})(u+v-1)/(u^2-v^2)\nonumber\\
&&- (DD-{\tilde {AA}})/(u+v),\label{228}
\end{eqnarray}
for the $U$-algebra. Here we use for brevity the following notations:
$[A,B]$ means the commutator $[A(u),B(v)]$, where the first parameter 
is $u$ and the second one is $v$; $DA$ stands for the noncommutative 
operator product $D(u)A(v)$; and ${\tilde {DA}}$ signifies $D(v)A(u)$
(where $v$ is the first parameter), and so on.   

\begin{theorem}\label{t6} 
Let $W$ be a complex vector space on which the QISM I algebra
acts by an algebra representation. Suppose ${\cal D}$ is a subset
of $\CC$ of the form $\{u_0+m \mid m\in\ZZ,\allowbreak j_-<m<j_+\}$,
where $u_0\in\CC$ and $j_\pm=\pm\iy$ or integer, such that

(i) $\{w\in W\mid C(u)w=0\}$ is 1-dimensional for any $u\in {\cal D}$,

(ii) if $u\in {\cal D}$, $0\neq w\in W$ and $C(u)w=0$ then
\begin{equation}
A(u)w\quad \left\{\quad \matrix{
\neq 0,& u\neq u_0+j_-+1,\cr
 =0, & u= u_0+j_-+1,}\right.
\label{2x1}\end{equation}
\begin{equation}
D(u)w\quad \left\{\quad \matrix{
\neq 0, & u\neq u_0+j_+-1,\cr
 =0,& u= u_0+j_+-1.}\right.
\label{2x2}\end{equation}
For each $u\in {\cal D}$ choose $0\neq F(u)\in W$ such that $C(u)F(u)=0$.
Then 
\begin{eqnarray}
{A(u)F(u)}&=&\Delta_-(u-\tfrac12)F(u-1),\qquad
u\in {\cal D},\quad u\neq u_0+j_-+1,\label{2x}\\
{D(u)F(u)}&=&\Delta_+(u+\tfrac12)F(u+1),\qquad
u\in {\cal D},\quad u\neq u_0+j_+-1,
\label{2xx}\end{eqnarray}
for certain scalar functions $\Delta_\pm (u\pm\tfrac12)$.
Furthermore, the operator
$\Delta(u)$, when acting on ${\rm Span}\,\{F(v)\mid  v\in {\cal D}\}$,
is scalar for $u\in\FSD\pm\tfrac12$ and it satisfies
\begin{equation}
\Delta (u)=\quad \left\{\quad \matrix{
\Delta_+(u)\Delta_-(u)\,,& u\in\FSD\pm\tfrac12\,,&
u\neq u_0+j_\pm\mp\tfrac12\,,\cr
0\,,&\qquad u=u_0+j_\pm\mp\tfrac12\,.&}\right.
\label{2x3}\end{equation}
\end{theorem}

\Proof
The matrix $L(u)$
satisfies all relations (\ref{21})--(\ref{213}). 
By substitution of $u=v-1$ in (\ref{25}) we get
$$C(v-1)A(v)=A(v-1)C(v).$$ 
Apply both sides to $F(v)$ ($v\in {\cal D}$, $v\neq u_0+j_0+1$), 
then we get the equation
$$C(v-1)A(v)F(v)=0.$$
$F(v-1)$ spans the zero space of $C(v-1)$ in $W$.
Thus
\begin{equation}A(v)F(v)\sim F(v-1).\label{229}\end{equation}
When we handle the commutator (\ref{28}) in a similar way we get 
$D(u)F(u)\sim F(u+1)$.
We write the proportionality factors as in (\ref{2x})--(\ref{2xx}),
by scalar factors $\Delta_\pm$ depending on $u$. 
Now apply the quantum determinant $\Delta(u-\tfrac12)$,
expressed by the second
formula of (\ref{9}), to $F(u)$ ($u\in {\cal D}$) and use (\ref{2x1})
or (\ref{2x2}) or (\ref{2x})--(\ref{2xx}). 
This yields (\ref{2x3}), with $u$ replaced by $u-\tfrac12$, and with
both sides acting on $F(u)$. Since $\Delta(u-\tfrac12)$ commutes 
with $A(v)$ and $D(v)$, the general case of (\ref{2x3}) then follows.
\eof

\begin{theorem}\label{t7}
Let $W$ be a complex vector space on which the QISM II algebra
acts by an algebra representation. Keep the other assumptions
of Theorem 1, except that $D(u)$ in (\ref{2x2}) is replaced by 
$-A(-u)$. Then the conclusions of Theorem 1 remain valid,
except that $D(u)$ in (\ref{2xx}) is replaced by $-A(-u)$.
\end{theorem}
\Proof
Analogous to the proof of Theorem \ref{t6}. Equation (\ref{229}) is now
obtained from (\ref{220}), while we get from (\ref{223}) the equation
$C(u+1)(2uD(u)- A(u))F(u)=0$.
In view of (\ref{141}) this implies $-A(-u)F(u)\sim F(u+1)$.
Use the second formula of (\ref{10}) for the proof of (\ref{2x3}).
\eof

\section{Some formulas for the classical special functions}
\setcounter{equation}{0}
For special functions of hypergeometric type there exists a large number
of formulas in which a (usually first order) differential operator
acting on the special function yields a special function of similar type
but with some parameters shifted. Usually such formulas occur in pairs,
with shifting of parameters in opposite directions. Below we list some
pairs of shift operator actions for which we will later give interpretations
in the framework of QISM I or II algebras. Throughout we use $u$ for the
parameter which is shifted. We give the formulas for the case of infinite
power series. For terminating power series the formulas can be rewritten
in terms of Jacobi polynomials, etc.

\subsection{Gauss hypergeometric function, Legendre function and
Jacobi polynomials}
The {\em Gauss hypergeometric function} \cite[Ch.2]{ba55}
\begin{equation}
{}_2F_1(a,b;c;x)\equiv F(a,b;c;x)=
\sum_{k=0}^\infty{(a)_k(b)_k\over (c)_kk!}\,x^k,\quad
|x|<1,\;c\ne0,-1,-2,\ldots\;,
\label{t66}
\end{equation}
is, up to a constant factor, the only analytic solution $f(x)$ in a
neighbourhood of 0
of the equation
\begin{equation}
(x(1-x)\partial_x^2+(c-(a+b+1)x)\partial_x-ab)\,f(x)=0.
\label{314}
\end{equation}
The solution is normalized by $f(0)=1$.
The function (\ref{t66}) has a unique analytic continuation to
$\CC\backslash[1,\iy)$. A second solution to (\ref{314}) is given by
\begin{equation}
f(x)=F(a,b;a+b-c+1;1-x),\quad c-a-b\ne1,2,\ldots\;.
\label{t100}
\end{equation}
Note also
\begin{equation}
F(a,b;c;x)=(1-x)^{-a}\,F(a,c-b;c;x/(x-1)).
\label{t101}
\end{equation}

{\em Jacobi polynomials}:
\begin{eqnarray}
P_n^{(\al,\be)}(x)=
{(\al+1)_n\over n!}\, F(-n,n+\al+\be+1;\al+1;\tfrac12(1-x))
=(-1)^n\,P_n^{(\be,\al)}(-x),
\label{t67}\\
\qquad\qquad\qquad\qquad
n=0,1,2,\ldots,\quad \al,\be\ne-1,-2,\ldots\;.\nonumber
\end{eqnarray}

{\em Legendre function} \cite[Ch.3]{ba55}:
\begin{eqnarray}
P_\nu^\mu(x)={2^\mu(x^2-1)^{-\mu/2}\over \Gamma(1-\mu)}\,
F(1-\mu+\nu,-\mu-\nu;1-\mu;\tfrac12-\tfrac12x),\nonumber\\
\qquad\qquad\qquad
x\in\CC\backslash(-\iy,1],\quad\mu\ne1,2,\ldots,\label{t98}\\
P_\nu^\mu(x)={\Ga(\nu+\mu+1)\over\Ga(\nu-\mu+1)}\,P_\nu^{-\mu}(x),\quad
\mu=1,2,\ldots\;.\label{t91}
\end{eqnarray}

Shift operator pairs:
\begin{eqnarray}
(x(1-x)\partial_x-bx+c-a-u)\,F(a+u,b;c;x)\qquad\qquad\nonumber\\
=(c-a-u)\,F(a+u-1,b;c;x),
\label{315}\\
(x\partial_x+a+u)\,F(a+u,b;c;x)=(a+u)\,F(a+u+1,b;c;x);\label{316}\\
\noalign{\medskip}
\Bigl(x\pa_x-{{b-c}\over{1-x}}+b-1+u\Bigr)\,\left[(1-x)^{a+u}
F(a+u,b+u;c+u;x)\right]\qquad\qquad\nonumber\\
=(c+u-1)\,(1-x)^{a+u-1}F(a+u-1,b+u-1;c+u-1;x),
\label{vv1}\\
((1-x)\pa_x+a+u)\,\left[(1-x)^{a+u}
F(a+u,b+u;c+u;x)\right]\qquad\qquad\qquad\nonumber\\
=\frac{(a+u)(b+u)}{c+u}\,(1-x)^{a+u+1}F(a+u+1,b+u+1;c+u+1;x);
\label{vv2}
\end{eqnarray}

\begin{eqnarray}
{\left(x(1-x)\pa_x+(\tfrac12-a)x+\tfrac12c-\tfrac12+(u-\tfrac12)(x-\tfrac12)+
{\delta\over u-\tfrac12}\right)\,F(a+u,a-u;c;x)}\nonumber\\
={(a-u)(a-c+u)\over 2(u-\tfrac12)}\,F(a+u-1,a-u+1;c;x),\qquad\label{30001}\\
{\left(-x(1-x)\pa_x-(\tfrac12-a)x-\tfrac12c+\tfrac12+(u+\tfrac12)(x-\tfrac12)
+{\delta\over u+\tfrac12}\right)\,F(a+u,a-u;c;x)}\nonumber\\
={(a+u)(a-c-u)\over 2(u+\tfrac12)}\,F(a+u+1,a-u-1;c;x),\qquad\label{30002}\\
\noalign{\hbox{where}}
\delta=\tfrac12(a-c+\tfrac12)(a-\tfrac12);\qquad\qquad\qquad\qquad
\label{t92}
\end{eqnarray}
\begin{eqnarray}
{\Bigl((x^2-1)^{1/2}\,\partial_x+{ux\over (x^2-1)^{1/2}}
\Bigr)P_\nu^u(x)}&=&(\nu+u)(\nu-u+1)P_\nu^{u-1}(x),
\label{v2}\\
{\Bigl((x^2-1)^{1/2}\,\partial_x-{ux\over (x^2-1)^{1/2}}
\Bigr)P_\nu^u(x)}&=&P_\nu^{u+1}(x).
\label{v3}
\end{eqnarray}

The following special pair of shift operator actions can be derived from
(\ref{vv1})--(\ref{vv2}). Fix $n\in\{1,2,\ldots\}$.
Define functions $F_k$ ($k\in\ZZ$) by
\begin{eqnarray}
F_k(x)={(n+k)!\over k!}\,(1-x)^{-n+k}\,F(-n+k,n+k+1;k+1;x),\quad
k=0,1,2,\ldots,\label{t93}\\
F_{-k}(x)=(-1)^k\,x^k\,(1-x)^{-k}\,F_k(x),\quad
k=1,2,\ldots\;.\quad\label{t94}
\end{eqnarray}
Then
\begin{eqnarray}
\Bigl(x\pa_x-{n\over 1-x}+n+k\Bigr)\,F_k(x)&=&(n+k)\,F_{k-1}(x),\label{t95}\\
((1-x)\pa_x-n+k)\,F_k(x)&=&-(n-k)\,F_{k+1}(x).\label{t96}
\end{eqnarray}

\subsection{Confluent hypergeometric function and Laguerre polynomials}
The {\it confluent hypergeometric function} \cite[\S6.3]{ba55}
\begin{equation}
{}_1F_1(a;c;x)\equiv
\Phi(a,c;x)=\sum_{k=0}^\infty {(a)_k\over (c)_k k!}\,x^k,\quad
c\ne0,-1,-2,\ldots\;,
\label{t69}
\end{equation}
is, up to a constant factor, the only entire analytic solution $f(x)$
of the equation
\begin{equation}
(x\partial_x^2+(c-x)\partial_x-a)\,f(x)=0.
\label{31}
\end{equation}
The solution is normalized by $f(0)=1$.
The confluent hypergeometric
function can be obtained as a limit case of the Gauss
hypergeometric function:
\begin{equation}
{}_1F_1(a;c;x)=\lim_{b\to\iy}{}_2F_1(a,b;c;b^{-1}x).
\label{t99}
\end{equation}
The other special functions we will discuss can also be obtained as limits
of the ${}_2F_1$-function. Accordingly, all further shift operator pairs
listed below are limit cases of shift operator pairs in \S3.1.

For arbitrary $a,c$
{\it Tricomi's $\Psi$-function} $\Psi(a,c;x)$
can be defined, for instance, by the contour integral representation
\cite[6.11(9)]{ba55}. It is, up to a constant factor,
the only analytic solution $f(x)$ of the equation (\ref{31}) on $(0,\iy)$
such that $f(x)$ is of at most polynomial growth as $x\to\iy$.
This characterization can be extracted from
\cite[\S6.7 and \S6.13.1]{ba55}.
All derivatives of this function $f$ have the same growth property as $f$.
The solution is normalized by
$f(x)=x^{-a}+\FSO(|x|^{-a-1})$ as $x\to\iy$.

{\it Laguerre polynomials}:
\begin{eqnarray}
L_n^\alpha (x)=
{(\alpha+1)_n\over n!}\,\Phi(-n,\alpha+1;x)=
{(-1)^n\over n!}\,\Psi(-n,\alpha+1;x),\label{t26}\\
n=0,1,2,\ldots\;.\nonumber
\end{eqnarray}

Shift operator pairs:
\begin{eqnarray}
{(x\partial_x-x+c-a-u)\,\Phi(a+u,c;x)}=(c-a-u)\,\Phi(a+u-1,c;x),
\label{32}\\
{(x\partial_x+a+u)\,\Phi(a+u,c;x)}=(a+u)\,\Phi(a+u+1,c;x);\quad\label{33}\\
\noalign{\medskip}
(x\partial_x-x+c+u-1)\,\Psi(u,c+u;x)=-\Psi(u-1,c+u-1;x),\quad\label{351}\\
\partial_x\,\Psi(u,c+u;x)=-u\,\Psi(u+1,c+u+1;x);\quad\label{352}
\end{eqnarray}

\newpage
\begin{eqnarray}
\left(\partial_x+\frac{u}{x}+{\de\over u-\tfrac12}\right)\left[
x^ue^{-\tfrac{x}{2}}\Psi(2\de+\tfrac12+u,2u+1;x)\right]\qquad\qquad
\nonumber\\
={2\de+{1\over 2}-u\over 2u-1}\,x^{u-1}e^{-\tfrac{x}{2}}
\Psi(2\de+\tfrac12+u-1,2u-1;x),\label{last3}\\
\left(-\partial_x+\frac{u}{x}+{\de\over u+\tfrac12}\right)\left[
x^ue^{-\tfrac{x}{2}}\Psi(2\de+\tfrac12+u,2u+1;x)\right]\qquad\qquad
\nonumber\\
={2\de+{1\over 2}+u\over 2u+1}\,x^{u+1}e^{-\tfrac{x}{2}}
\Psi(2\de+\tfrac12+u+1,2u+3;x).\label{last4}
\end{eqnarray}

\subsection{Parabolic cylinder function and Hermite polynomials}
The {\it parabolic cylinder function} \cite[\S8.2]{ba55} can be defined by
\begin{equation}
D_\nu(x)=2^{\frac12(\nu-1)}e^{-x^2/4}x\,
\Psi(\tfrac12-\tfrac12\nu,\tfrac32;\tfrac12x^2).
\label{t70}
\end{equation}
The function $f(x)=e^{\frac14x^2}D_\nu(x)$
is, up to a constant factor,
the only entire analytic solution of the equation
\begin{equation}
(\pa_x^2-x\pa_x+\nu)\,f(x)=0
\label{t71}
\end{equation}
such that $f(x)$ is on $(0,\iy)$ of at most polynomial growth as
$x\to\iy$.
All derivatives of this function $f$ have the same growth property as $f$.
The solution is normalized by $f(x)=x^\nu+\FSO(x^{\nu-2})$ as $x\to\iy$.

{\it Hermite polynomials} (polynomials of degree $n$):
\begin{equation}
H_n(x)=2^{\frac12n}\,e^{\frac12 x^2}\,D_n(2^{\frac12}x),\quad
n=0,1,2,\ldots\;.
\label{t72}
\end{equation}

Shift operator pairs:
\begin{eqnarray}
(\pa_x-x)\,(e^{\frac14x^2}D_{-u}(x))&=&-e^{\frac14x^2}D_{-u+1}(x),
\label{t63}\\
\pa_x(e^{\frac14x^2}D_{-u}(x))&=&-u\,e^{\frac14x^2}D_{-u-1}(x).\label{t62}
\end{eqnarray}

\subsection{Bessel function}
The {\it Bessel function} \cite[Ch.7]{ba55}
\begin{eqnarray}
J_\nu(x)=\sum_{m=0}^\infty{(-1)^m(\tfrac12x)^{2m+\nu}
\over m!\,\Gamma(m+\nu+1)}={(\tfrac12x)^\nu\over \Gamma
(\nu+1)}\,{}_0F_1(-;\nu+1;-\tfrac14 x^2),\qquad\nonumber\\
\nu\ne-1,-2,\ldots,\label{350}\\
J_\nu(x)=(-1)^\nu\,J_{-\nu}(x),\quad\nu=-1,-2,\ldots,\label{t97}
\end{eqnarray}
is a solution $f(x)$ of the equation
\begin{equation}
((x\partial_x)^2+x^2-\nu^2)\,f(x)=0.
\label{353}
\end{equation}
Hence the function $f(x)={}_0F_1(-;1+\nu;\tfrac14 x^2)$\quad
($\nu\ne-1,-2,\ldots$) satisfies the equation
\begin{equation}
(\partial_x^2+(2\nu+1)x^{-1}\partial_x-1)\,f(x)=0.\label{355}
\end{equation}
Another solution to (\ref{355}) (for any $\nu$) is given by
$f(x)=x^{-\nu}\,K_\nu(x)$ ($x>0$), where $K_\nu$ is the
{\it modified Bessel function of the third kind}
defined for instance by the integral representation
\cite[7.12(23)]{ba55}. It is, up to a constant factor, the unique 
analytic solution
$f(x)$ of (\ref{355}) on $(0,\iy)$
which tends to 0 faster than any inverse power as
$x\to\iy$ (cf.\ \cite[7.13(7)]{ba55}).
All derivatives of this function $f$ have the same growth property as $f$.
The solution is normalized by
$f(x)=(\pi/2)^{\frac12}x^{-\nu-\frac12}e^{-x}(1+\FSO(x^{-1}))$ as $x\to\iy$.

Shift operator pairs:
\begin{eqnarray}
(x\partial_x+2u)\,[x^{-u}K_u(x)]&=&
-x^{-u+1}K_{u-1}(x),\label{t65}\\
x^{-1}\partial_x\,[x^{-u}K_u(x)]&=&-x^{-u-1}K_{u+1}(x);
\label{356}\\
\noalign{\medskip}
(\partial_x+u x^{-1})J_u(x)&=&
J_{u-1}(x),\label{354}\\
(-\partial_x+u x^{-1})J_u(x)&=&
J_{u+1}(x).\label{t64}
\end{eqnarray}

\section{Lie algebras of first order differential operators}
\setcounter{equation}{0}
In the QISM I case  we obtain from (\ref{211}) that
\begin{eqnarray}
D(u-1)A(u)-A(u)D(u-1)&=&B(u-1)C(u)-B(u)C(u-1),\nonumber\\
D(u)A(u+1)-A(u+1)D(u)&=&B(u)C(u+1)-B(u+1)C(u).\nonumber
\end{eqnarray}
In combination with (\ref{9}) (fourth resp.\ second formula) this yields
\begin{eqnarray}
D(u-1)A(u)&=&B(u-1)C(u)+\De(u-\tfrac12),\label{t32}\\
A(u+1)D(u)&=&B(u+1)C(u)+\De(u+\tfrac12).\label{t33}
\end{eqnarray}
Assume now that we have a representation of the QISM I algebra on a space
of functions in one variable, analytic on a certain region, such that
(i) $B(u)=B_0$ is independent of $u$,
(ii) $\De(u)$ is scalar for all $u$,
(iii) $A(u)$ and $D(u)$ are first order differential operators.
Then the equations (\ref{t32})--(\ref{t33}) show that the second order
operator $B_0C(u)$ has a suitable form for the factorization method,
which was originated by Schr\"odinger and was due in its definitive form to
Infeld and Hull \cite{infeld}.

Similarly, in the QISM II case we obtain from (\ref{215}) that
\begin{eqnarray}
A(u+1)A(-u)-A(-u)A(u+1)&=&-B(u+1)C(-u)+B(-u)C(u+1),\nonumber\\
A(u)A(-u+1)-A(-u+1)A(u)&=&-B(u)C(-u+1)+B(-u+1)C(u).\nonumber
\end{eqnarray}
Assume that $B(-u)=B(u)$ and $C(-u)=C(u)$ (part of the symmetry properties
(\ref{141})).
Then it follows in combination with (\ref{10}) (fourth resp.\ second formula)
that
\begin{eqnarray}
-A(-u+1)A(u)&=&B(u-1)C(u)+\De(u-\tfrac12),\label{t34}\\
-A(u+1)A(-u)&=&B(u+1)C(u)+\De(u+\tfrac12).\label{t35}
\end{eqnarray}
So assume that we have a representation of the QISM II algebra on a space
of functions in one variable, analytic on a certain region, such that
(i) $B(u)=B_0$ is independent of $u$,
(ii) $C(u)=C(-u)$ for all $u$,
(iii) $\De(u)$ is scalar for all $u$,
(iv) $A(u)$ is a first order differential operator.
Then equations (\ref{t34})--(\ref{t35}) show that the second order operator
$B_0C(u)$ has a suitable form for the Infeld-Hull factorization method.

The factorization method is summarized in Miller \cite[Ch.\ 7]{mi68}.
Under some special assumptions on the type of factorizing operators
(first order part not depending on $u$, zero order part of degree at most
one in $u$), a complete classification of all possibilities is given.

These factorizing operators of special type give rise to a Lie algebra of
first order differential operators in two variables. Indeed,
assume that we have elements
$A(u), D(u), K(u), \De(u)$ ($u\in\CC$) of an associative algebra $\FSA$
such that
$\De(u)$ is scalar and
\begin{eqnarray}
D(u-1)A(u)&=&K(u)+\De(u-\tfrac12),\label{t36}\\
A(u+1)D(u)&=&K(u)+\De(u+\tfrac12)\label{t37}.
\end{eqnarray}
Assume that $A(u)$ and $D(u)$ are of degree at most one in $u$:
\begin{equation}
A(u)=A_0+A_1u,\quad
D(u)=D_0+D_1u.
\label{t38}
\end{equation}
Consider now the algebra spanned by elements of the form
$B t^k\pa_t^l$, where $B\in\FSA$. Define in this
algebra the elements
\begin{equation}
J^+=t(D_0+D_1\,t\pa_t)=tD(t\pa_t),\qquad
J^-=t^{-1}(A_0+A_1\,t\pa_t)=t^{-1}A(t\pa_t).
\label{t61}
\end{equation}
Then, in view of (\ref{t36})--(\ref{t37}), we have
$$
J^+J^-=K(t\pa_t)+\De(t\pa_t-\tfrac12),\qquad
J^-J^+=K(t\pa_t)+\De(t\pa_t+\tfrac12).
$$
Hence
\begin{equation}
[J^-,J^+]=\De(t\pa_t+\tfrac12)-\De(t\pa_t-\tfrac12).
\label{t39}
\end{equation}
Assume furthermore that $A_1$ and $D_1$ commute.
Then it follows from (\ref{t38}) and (\ref{t39}) that, for certain
complex constants $Q_2$ and $Q_1$, we have
\begin{equation}
[J^-,J^+]=2Q_2t\pa_t+Q_1.
\label{t40}
\end{equation}
Clearly, we have also the commutators
\begin{equation}
[t\pa_t,J^\pm]=\pm J^\pm.
\label{t41}
\end{equation}
Thus the elements $J^+,J^-,t\pa_t$ span, together with the central element $1$,
a four-di\-men\-sional complex Lie algebra denoted $\FSG(a,b)$
in Miller \cite[\S 2-5]{mi68}. Here $a^2=-Q_2$ and $b=Q_1$.
They fall apart into three isomorphism classes:
(i) sl$(2,\CC)\oplus\CC$ ($a\ne0$), (ii) complexification of the
harmonic oscillator algebra, i.e.\ of the semidirect sum of $\RR$ with
the Heisenberg Lie algebra
($a=0$, $b\ne0$), (iii) e$(2,\CC)\oplus\CC$ ($a=b=0$),
where e$(2,\CC)$
is the complexified Lie algebra of the group of plane motions.

Conversely, if $J^\pm$ are of the form (\ref{t61}) and if (\ref{t40}) holds
then
$$
A(u+1)D(u)-D(u-1)A(u)=2Q_2u+Q_1.
$$
If we then put
$$
\De(u)=Q_2u^2+Q_1u+Q_0
$$
for some constant $Q_0$ then $\De(u+\tfrac12)-\De(u-\tfrac12)=2Q_2u+Q_1$
and
$$
A(u+1)D(u)-\De(u+\tfrac12)=D(u-1)A(u)-\De(u-\tfrac12).
$$
So, if $K(u)$ is put equal to the left hand side of the above identity,
we recover (\ref{t36})--(\ref{t37}).

Let us return to equations (\ref{t36})--(\ref{t37}).
Suppose that the algebra $\FSA$ acts on some linear space $W$ and that,
for some $u_0\in\CC$ and some $F(u_0)\in W\backslash\{0\}$, we have that
$K(u_0)\,F(u_0)=0$. Then it follows from (\ref{t36})--(\ref{t37}) that
$K(u_0-1)\,(A(u_0)F(u_0))=0$ and $K(u_0+1)\,(D(u_0)F(u_0))=0$.
Moreover, $D(u_0-1)A(u_0)F(u_0)=\De(u_0-\tfrac12)\,F(u_0)$ and
$A(u_0+1)D(u_0)F(u_0)=\De(u_0+\tfrac12)\,F(u_0)$.
Thus, for $k=1,2,\ldots$,
we can recursively define $F(u_0+k)=\const D(u_0+k-1)F(u_0+k-1)$ and
$F(u_0-k)=\const A(u_0-k+1)F(u_0-k+1)$, as long as these vectors are nonzero.
In that way we obtain strings of vectors $F(u)$ ($u\in \FSD$) as in Theorem
\ref{t6}, where $K(u)\,F(u)=0$ for $u\in\FSD$ and
(\ref{2x})--(\ref{2xx}) are valid for a certain choice of $\De_\pm(u)$.
If, moreover, the operators $J^\pm$ are defined by (\ref{t61}) then
$$
J^\pm(t^uF(u))=\De_\pm(u\pm\tfrac12)\,t^{u\pm1}F(u\pm1),\quad u\in\FSD.
$$
So we have a Lie algebra acting on the elements $t^u F(u)$.
The algebra acting on the elements $F(u)$ is not a Lie algebra, in general,
but we will see later in this paper that this algebra action can often be
extended to a QISM I algebra action. The crucial point in making this
extension is to find an operator $C_0$ acting on $W$ such that
$C_0 F(u)=-uF(u)$ for $u\in\FSD$. A necessary condition for finding such an
operator will be that the elements $F(u)$ ($u\in\FSD$) are
linearly independent.
Note that this is not yet guaranteed by the above assumptions.
But, of course, we do know that the elements $t^u F(u)$ ($u\in\FSD$)
are linearly independent.

In \cite[\S 2-7]{mi68} Miller assumes that
\begin{eqnarray}
J^+&=&t(\pa_x-k(x)t\pa_t+j(x)),\label{t42}\\
J^-&=&t^{-1}(-\pa_x-k(x)t\pa_t+j(x))\label{t43}
\end{eqnarray}
for certain analytic functions $k$ and $j$.
Then he shows that (\ref{t40}) holds if and only if
\begin{eqnarray}
&&k'(x)+k(x)^2=Q_2,\label{t44}\\
&&j'(x)+k(x)j(x)=-\tfrac12 Q_1.\label{t45}
\end{eqnarray}
The general solution of  equations (\ref{t44})--(\ref{t45}) will depend
on two parameters, but one parameter is trivial because the equations are
invariant under translation.
Solution of the equations yields six different cases,
two for each isomorphism class of the Lie algebra $\FSG(a,b)$, depending
on whether $k(x)$ is constant or not.
The list of solutions is as follows (cf.\ \cite[p.272]{mi68}, we do not
give the trivial translation parameter):
\begin{eqnarray}
&{\rm Type}\;(A)\quad&
k(x)=a\cot ax,\quad
j(x)={b\over 2a}\,\cot ax+{q\over \sin ax}\,.\label{t46}\\
&{\rm Type}\;(B)\quad&
k(x)=ia,\quad
j(x)={ib\over 2a}+qe^{-iax}.\label{t47}\\
&{\rm Type}\;(C')\quad&
k(x)=x^{-1},\quad
j(x)=-\tfrac14bx+qx^{-1}.\label{t48}\\
&{\rm Type}\;(D')\quad&
k(x)=0,\quad
j(x)=-\tfrac12bx.\label{t49}\\
&{\rm Type}\;(C'')\quad&
k(x)=x^{-1},\quad
j(x)=qx^{-1}.\label{t50}\\
&{\rm Type}\;(D'')\quad&
k(x)=0,\quad
j(x)=q.\label{t51}
\end{eqnarray}
Here $a,b$ are the parameters from $\FSG(a,b)$ and $q$ is another parameter.
For types $A$ and $B$ we have $a\ne0$, for types $C'$ and $D'$ we have
$a=0$ and $b\ne0$, and for types $C''$ and $D''$ we have $a=b=0$.
For types $A$, $C'$ and $C''$ $k(x)$ is not constant, but for
types $B$, $D'$ and $D''$ it is. In the following we do not consider
the trivial case $D''$ because it does not give any shift operator pair.

The operators (\ref{t42})--(\ref{t43}) are in a certain normal form.
We want to transform them into another normal form which is better adapted
to the QISM I algebra. This is done in the following lemmas, which can be
proved by straightforward computation.

\begin{lemma}\label{t52}
Let $J^\pm$ be given by (\ref{t42})--(\ref{t43}) and assume that
(\ref{t40}), and thus (\ref{t44})--(\ref{t45}), hold.
Make a transformation of the variables $x,t$ by replacing $t$ by
$(\psi(x))^{-1} t$,
where $\psi$ is such that
$$
(\psi'(x)/\psi(x))^2=-Q_2+k(x)^2.
$$
Then
\begin{eqnarray}
J^+&=&t(D_{01}(x)\pa_x+D_{00}(x)+\de t\pa_t),\label{t53}\\
J^-&=&t^{-1}(A_{01}(x)\pa_x+A_{00}(x)+\al t\pa_t),\label{t54}
\end{eqnarray}
where $A_{00},A_{01},D_{00},D_{01}$ are certain analytic functions
with $A_{01}$ and $D_{01}$ not identically zero, and
$\al,\de$ are complex constants such that
$\al\de=Q_2$.

Furthermore, $A_{01}(x)/D_{01}(x)$ is constant or not depending on whether
$k(x)$ is constant or not.
\end{lemma}

\begin{lemma}\label{t60}
Let $J^\pm$ be given by (\ref{t53})--(\ref{t54}).
Write
\begin{equation}
A_0=A_{01}(x)\pa_x+A_{00}(x),\quad
D_0=D_{01}(x)\pa_x+D_{00}(x).
\label{t56}
\end{equation}
Then (\ref{t40}) holds if and only if $\al\de=Q_2$ and
\begin{equation}
[D_0,A_0]=\de A_0+\al D_0-Q_1.
\label{t55}
\end{equation}
Furthermore, (\ref{t55}) holds if and only if
\begin{eqnarray}
D_{01}(x)A_{01}'(x)-A_{01}(x)D_{01}'(x)&=&\de A_{01}(x)+\al D_{01}(x),
\label{t57}\\
D_{01}(x)A_{00}'(x)-A_{01}(x)D_{00}'(x)&=&\de A_{00}(x)+\al D_{00}(x)-Q_1.
\label{t58}
\end{eqnarray}
\end{lemma}

\begin{lemma}\label{t59}
Let $J^\pm$ be given by (\ref{t53})--(\ref{t54})
with $A_{01}$ and $D_{01}$ not identically zero
and assume that (\ref{t40}),
and thus (\ref{t57})--(\ref{t58}), hold.
Make a transformation of the variables $x,t$ by replacing $t$ by
$(\phi(x))^{-1} t$,
where $\phi$ is such that
$$
2\phi(x)\phi'(x)={\de-\al\phi(x)^2\over A_{01}(x)}\,.
$$
Then
\begin{eqnarray}
J^+&=&t(\chi(x)\pa_x+j^+(x)-k(x)t\pa_t),\nonumber\\
J^-&=&t^{-1}(-\chi(x)\pa_x+j^-(x)-k(x)t\pa_t),\nonumber
\end{eqnarray}
for certain analytic functions $\chi$ (not identically zero),
$j^\pm$ and $k$. Furthermore $k(x)$ is constant or not according to
whether $A_{01}(x)/D_{01}(x)$ is constant or not.
Finally, for a suitable analytic function $f$, not identically zero,
and after a suitable transformation of the $x$-variable the operators
$f(x)^{-1}\,J^\pm\circ f(x)$ take the
form (\ref{t42})--(\ref{t43}).
\end{lemma}

Equations (\ref{t56}) and (\ref{t55}) will not change if
the operators $A_0$ and $D_0$ are replaced by operators
$f(x)^{-1}\,A_0\circ f(x)$ and $f(x)^{-1}\,D_0\circ f(x)$, respectively,
for a suitable analytic function $f$, not identically zero.
We call such transformations {\it gauge transformations}.
Neither do the equations change when we make an analytic transformation
of the $x$-variable. We will consider solutions to the equations
as equivalent if they can be obtained from each other by the two types
of transformations just described.

In the formulas below we list operators $A(u)=A_0+\al u$,
$D(u)=D_0+\de u$ such that $A_0$ and $D_0$ have the form (\ref{t56})
and such that they satisfy the equivalent conditions of Lemma \ref{t60}.
These formulas can be derived either from (\ref{t46})--(\ref{t50})
by use of the above lemmas, or by straightforward verification that the
conditions of Lemma \ref{t60} are satisfied.
\begin{eqnarray}
{\rm Type}\;(A)\quad&&
A(u)=x(1-x)\pa_x-bx+c-a-u,\quad
D(u)=x\pa_x+a+u,\nonumber\\
&&[D_0,A_0]=-x^2\pa_x-bx,\quad
Q_1=c-2a.\label{t73}\\
{\rm Type}\;(B)\quad&&
A(u)=x\pa_x-x+c-a-u,\quad
D(u)=x\pa_x+a+u,\nonumber\\
&&[D_0,A_0]=-x,\quad
Q_1=c-2a.\label{t74}\\
{\rm Type}\;(C')\quad&&
A(u)=x\pa_x-x+c-1+u,\quad
D(u)=\pa_x,\nonumber\\
&&[D_0,A_0]=\pa_x-1,\quad
Q_1=1.\label{t75}\\
{\rm Type}\;(D')\quad&&
A(u)=\pa_x-x,\quad
D(u)=\pa_x,\nonumber\\
&&[D_0,A_0]=1,\quad
Q_1=-1.\label{t76}\\
{\rm Type}\;(C'')\quad&&
A(u)=x\pa_x+2u,\quad
D(u)=x^{-1}\pa_x,\nonumber\\
&&[D_0,A_0]=2x^{-1}\pa_x,\quad
Q_1=0.\label{t77}
\end{eqnarray}
For all these  types we can
give functions $F(u)$ on which $A(u)$ and $D(u)$ act as shifting operators.
See equations
(\ref{315})--(\ref{316}),
(\ref{32})--(\ref{33}),
(\ref{351})--(\ref{352}),
(\ref{t63})--(\ref{t62}),
(\ref{t65})--(\ref{356}),
respectively.
These functions $F(u)$ are not uniquely determined. We might write down
similar formulas with another choice $F(u)$ for the solution of
the corresponding second order equation.

Next we discuss a form of the operators $J^\pm$
(cf.\ (\ref{t42})--(\ref{t43})) which we will meet
in the case of the QISM II algebra.
Put
\begin{eqnarray}
J^+=t(-A_{01}(x)\pa_x-A_{00}(x)+A_1(x)\,(t\pa_t+\tfrac12)),\label{t85}\\
J^-=t^{-1}(A_{01}(x)\pa_x+A_{00}(x)+A_1(x)\,(t\pa_t-\tfrac12)),\label{t86}
\end{eqnarray}
where $A_{00},A_{01},A_1$ are certain analytic functions with $A_{01}$
not identically zero.
Observe that  (\ref{t42})--(\ref{t43}) is of the form (\ref{t85})--(\ref{t86})
if and only if
$j(x)=0$ (then necessarily, by (\ref{t45}), $Q_1=0$).
This occurs non-trivially in (\ref{t46})--(\ref{t51})
($k(x)$ is not constant and $j(x)=0$) precisely for Types $A$ ((\ref{t46})
with $b=q=0$) and $C''$ ((\ref{t50}) with $q=0$).
On the other hand operators $J^\pm$ of the form (\ref{t85})--(\ref{t86})
can be brought in the form (\ref{t42})--(\ref{t43}) by suitable gauge
transformation and $x$-transformation.

Write $A_0=A_{01}(x)\pa_x+A_{00}(x)$ as before.
The following lemma can be proved by straightforward computation.

\begin{lemma}\label{t87}
Let $J^\pm$ be given by (\ref{t85})--(\ref{t86}).
Then (\ref{t40}) holds if and only if $Q_1=0$ and
\begin{equation}
[A_0,A_1]+A_1^2=Q_2.
\label{t88}
\end{equation}
Furthermore, (\ref{t88}) holds if and only if
\begin{equation}
A_{01}A_1'+A_1^2=Q_2.
\label{t89}
\end{equation}
\end{lemma}

\bigbreak
In the formulas below we list operators $A(u)=A_0+A_1 u$
with $A_0$ an analytic operator of the form (\ref{t56})
and $A_1$ an analytic function
such that $A_{00},A_{01},A_1$  satisfy the equivalent conditions of
Lemma \ref{t87}. 
\begin{eqnarray}
{\rm Type}\;(A)\quad&&
A(u)=x(1-x)\pa_x+(u-a)(x-\tfrac12),\nonumber\\
&&[A_1,A_0]=-x(1-x).\label{v73}\\
{\rm Type}\;(C'')\quad&&
A(u)=ux^{-1}+\pa_x,\nonumber\\
&&[A_1,A_0]=x^{-2}.\label{v75}
\end{eqnarray}
For these two types we can
give functions $F(u)$ on which $A(u)$ and $-A(-u)$ act as shifting operators.
See equations (\ref{30001})--(\ref{30002}) ($c=a+\tfrac12$) and
(\ref{354})--(\ref{t64}), respectively.

\section{Rank 1 $L$-operators for the QISM I algebra}
\setcounter{equation}{0}
The $T$-algebra is the algebra with the matrix elements of $T(u)$ as
generators and with the relation (\ref{3}).
We pass to a quotient algebra by adding the relation
$T(u)=u\,T(1)+(1-u)\,T(0)$.
In other words, we make the ansatz that $T(u)$ is of the form
\begin{equation}
L(u)=\left(\matrix{A_1 u+A_0&B_1u+B_0\cr C_1u+C_0&D_1u+D_0}
\right)=L_1 u + L_0.
\label{41}
\end{equation}
Substitute (\ref{41}) in relation (\ref{3}).
Then we get the algebra with $A_i,B_i,C_i,D_i,\,i=0,1$ as generators
and with relations
\begin{eqnarray}
{L_1^{(1)}L_1^{(2)}}&=&L_1^{(2)}L_1^{(1)},\quad 
L_1^{(1)}L_0^{(2)}=L_0^{(2)}L_1^{(1)},\label{42}\\
{[L_0^{(1)},L_0^{(2)}]}&=&- [P,L_1^{(1)}L_0^{(2)}].
\label{43}\end{eqnarray}
The relations (\ref{42}) imply that the entries of the $L_1$-matrix
are in the center of the algebra. Let us pass once more to a quotient
algebra by adding the relation
\begin{equation}
L_1=\left(\matrix{\alpha&\beta\cr \gamma&\delta}\right)
\label{44}
\end{equation}
for certain $\al,\be,\ga,\de\in\CC$.
A representation of the algebra with relations (\ref{42}) and
(\ref{43}) which has the
property that all elements in the center of the algebra are represented as
scalars, can also be viewed as a representation of the algebra with
relations (\ref{43}) and (\ref{44}) for a certain choice of $\al,\be,\ga,\de$.

\begin{remark}\label{t8}
Consider for a moment the more general situation of
a QISM I algebra relation (\ref{3}) with $n$-dimensional auxiliary
space, i.e., 
$R(u)=u+ P$, where $P\,(x\otimes y)=y\otimes x$ for
$x,y\in \CC^n$. Just as in the $2\times 2$ case we make the ansatz that
$T(u)=uM+L$, where $M$ and $L$ are $n\times n$ matrices with
a priori non-commuting matrix elements.
Then the matrix entries of $M$ are in the center of the algebra.
We add the relations $M_{ij}=\mu_{ij}$ for certain $\mu_{ij}\in\CC$.
Then we obtain the algebra with the $L_{ij}$ as generators and with relations
\begin{equation}
[L^{(1)},L^{(2)}]=- [P,\mu^{(1)}L^{(2)}],
\label{40002}
\end{equation} 
hence
\begin{equation}
[L_{ir},L_{js}]=- \mu_{jr}L_{is}+ \mu_{is}L_{jr}.
\label{40003}
\end{equation}
Alternatively, we may consider the linear space $\gog$ with
the $L_{ij}$ ($i,j=1,\ldots,n$) as basis vectors (so they are
linearly independent) and with antisymmetric
bilinear product defined by (\ref{40003}).
We claim that, for any choice of the $\mu_{ij}$, the space $\gog$ equipped
with this product becomes a Lie algebra, i.e., the product satisfies the
Jacobi identity. Indeed, we have
$$
[[L_{ir},L_{js}],L_{kt}]=
\mu_{jr}\mu_{ks}L_{it}-\mu_{jr}\mu_{it}L_{ks}
-\mu_{is}\mu_{kr}L_{jt}+\mu_{is}\mu_{jt}L_{kr}
$$
and two similar identities obtained by cyclic permutation
of the indices $(ijk)$ and $(rst)$. Addition of the three identities
yields 0 on the right hand side.

The fact that $\gog$ is a Lie algebra is equivalent to saying that
the $L_{ij}$ are linearly independent in the algebra with generators $L_{ij}$
and relations (\ref{40003}). Then this algebra is the universal enveloping
algebra of the Lie algebra $\gog$.
\end{remark}

Let us return to the case $n=2$.
The commutator
(\ref{43}) yields a Lie algebra spanned by $A_0,B_0,C_0,$ and $D_0$.
The componentwise form of the commutator (\ref{43}) is
as follows:
\begin{eqnarray}
{[A_0,B_0]}&=&\beta A_0-\alpha B_0,
\quad [A_0,C_0]=\alpha C_0-\gamma A_0,\nonumber\\
{[A_0,D_0]}&=&\beta C_0-\gamma B_0,
\quad [B_0,C_0]=\alpha D_0-\delta A_0,\nonumber\\
{[B_0,D_0]}&=&\beta D_0-\delta B_0,
\quad [C_0,D_0]=\delta C_0-\gamma D_0.
\label{45}\end{eqnarray}
The quantum determinant (\ref{9}) now has the form
\begin{eqnarray}
\Delta(u)&=&(\alpha\delta-\beta\gamma)u^2+Q_1u+Q_0,\label{46}\\
Q_1&=&\alpha D_0-\gamma B_0+\delta A_0-\beta C_0,\label{t1}\\
Q_0&=&A_0D_0-B_0C_0+\tfrac12\left(\alpha D_0+\gamma B_0-
\delta A_0-\beta C_0\right)-\tfrac14(\alpha\delta-\beta\gamma).
\label{t2}\end{eqnarray}
Here the right hand sides of (\ref{t1}) and (\ref{t2}) give operators
in the center of the algebra. We consider (\ref{t1}) and (\ref{t2}) as
added relations, for a certain choice of $Q_1,Q_2\in\CC$.
So a representation of the algebra with relations (\ref{42}) and 
(\ref{43}) which has the
property that all elements in the center of the algebra are represented as
scalars, can also be viewed as a representation of the algebra with
relations (\ref{43}), (\ref{44}), (\ref{t1}) and (\ref{t2})
for a certain choice of $\al,\be,\ga,\de,Q_1,Q_2$.

{}From now on we assume that 
\begin{equation}
\beta=0,\quad\gamma\neq 0.
\label{4--extra}\end{equation}

\begin{remark}\label{t102}
We will determine the type of Lie algebra given by
(\ref{45}) with (\ref{t1}), (\ref{4--extra}).
These commutators can be equivalently written as
\begin{eqnarray}
{[\,-\ga^{-1}C_0\,,
\,A_0-\al\ga^{-1}C_0\,]}&=&-(A_0-\al\ga^{-1}C_0),\nonumber\\
{[\,-\ga^{-1}C_0\,,
\,D_0-\de\ga^{-1}C_0\,]}&=&D_0-\de\ga^{-1}C_0,\nonumber\\
{[\,A_0-\al\ga^{-1}C_0\,,
\,D_0-\de\ga^{-1}C_0\,]}&=&2\al\de(-\ga^{-1}C_0)+Q_1.
\nonumber\end{eqnarray}
These three equations have the same structure as the equations
(\ref{t40}), (\ref{t41}), which we took from Miller's book
\cite{mi68} and which give rise to a Lie algebra $\FSG(a,b)$ with
$a^2=-\al\de$ and $b=Q_1$.
Thus we find the same three types of Lie algebras spanned by
$A_0,C_0,D_0$ and the central element 1 as in the discussion after
(\ref{t41}).

In the following we will obtain realizations of these Lie algebras
as operators acting on functions of one variable. Here $A_0$ and $D_0$
will be first order differential operators, but $C_0$ a second order
differential operator or an integro-differential operator.
It is interesting to compare this with Miller \cite{mi68},
whose only realizations of these Lie algebras by operators acting on functions
of one variable are by first order differential operators.
\end{remark}

The following lemma can be proved in a straightforward way. It
shows that equations (\ref{45}), (\ref{t1}) and (\ref{t2}), with
(\ref{4--extra}), and under the assumption that $B_0$ is injective,
can be equivalently written in a much more simple form.

\begin{lemma}\label{t9}
Let $\alpha,\gamma,\delta,Q_0,Q_1$ be scalars, with $\gamma\neq 0$.
Let $A_0,B_0,C_0,D_0$
be operators acting on some linear space $W$. Let $B_0$ be injective.
Then the following three statements are equivalent:
\vskip 2mm

(a)\quad
$\pmatrix{A_0+\al u&B_0\cr C_0+\ga u&D_0+\de u}$ 
\vskip 2mm

is a representation of the QISM I algebra with quantum determinant 

$\De(u)=\al\de u^2+Q_1 u+Q_0$;

(b)\quad
The six commutators (\ref{45}) and formulas (\ref{t1}), (\ref{t2})
are valid with $\be=0$;

(c)\quad
The following three equalities are valid:
\begin{eqnarray}
[D_0,A_0]&=&\alpha D_0+\delta A_0-Q_1,\label{47}\\
\ga B_0&=&\alpha D_0+\delta A_0-Q_1,\label{4b0}\\
B_0C_0&=&(A_0+\al)D_0-(\tfrac14\alpha\delta+\tfrac12Q_1+Q_0).\label{t3}
\end{eqnarray}
Moreover,
if $\{A_0,B_0,C_0,D_0,\al,\ga,\de,Q_0,Q_1\}$ satisfy these equivalent
conditions, then so do
$\{\la\mu A_0,\la\mu\nu^{-1}B_0,\la\nu C_0,\la D_0,\la\mu\al,
\la\nu\ga,\la\de,\la^2\mu Q_0,\la^2\mu Q_1\}$ (where $\la,\mu,\nu$
are nonzero scalars) and
$\{D_0,B_0,C_0,A_0,-\de,-\ga,-\al,Q_0,-Q_1\}$.
\end{lemma}

The next proposition is, in a certain sense, an inverse to Theorem \ref{t6}.
If operators $A(u)=A_0+\al u$
and $D(u)=D_0+\de u$ act on basis vectors $F(u)$ as in
(\ref{2x})--(\ref{2xx}) (part of the conclusion of Theorem \ref{t6})
then we can define actions of operators $B_0$ and $C_0$ such that
condition (c) of Lemma \ref{t9} is satisfied for certain $Q_0$ and $Q_1$.
So we then have obtained a representation of the QISM I algebra.

\begin{prop}\label{t90}
Let $\FSD$ be a subset of $\CC$ of the same form as in Theorem \ref{t6}.
Let $W$ be a complex vector space spanned by linearly independent vectors
$F(u)$ ($u\in\FSD$). Let $\al,\de$ be scalars.
Let $A_0$ and $D_0$ be linear operators on $W$ such that
$A(u)=A_0+\al u$ and $D(u)=D_0+\de u$ act on $F(u)$ as in
(\ref{2x})--(\ref{2xx}). Let $\De(u)$ ($u\in\FSD\pm\tfrac12$) be defined
by (\ref{2x3}) and assume that it has the form
$$\De(u)=\al\de u^2+Q_1u+Q_0$$
for certain scalars $Q_0,Q_1$.
Then (\ref{47}) is valid. Now define $B_0$ by (\ref{4b0}) (with $\ga=1$)
and $C_0$ by
$C_0F(u)=-uF(u)$ ($u\in\FSD$). Then (\ref{t3}) is satisfied
and also the other commutators in (\ref{45}) (with $\be=0$, $\ga=1$).
Then condition (a) of Lemma \ref{t9} is also satisfied.
(However, $B_0$ is not necessarily injective.)
\end{prop}

\Proof
For $u\in\FSD$ we find
\begin{eqnarray}
A_0D_0\,F(u)&=&\left(
u(Q_1-\al D_0-\de A_0)-\al D_0+(\tfrac14\al\de+\tfrac12Q_1+Q_0)
\right)F(u),\nonumber\\
D_0A_0\,F(u)&=&\left(
u(Q_1-\al D_0-\de A_0)+\de A_0+(\tfrac14\al\de-\tfrac12Q_1+Q_0)
\right)F(u).\nonumber
\end{eqnarray}
Hence (\ref{47}) is satisfied when both sides act on $F(u)$.
For $B_0$ and $C_0$ as defined in the proposition, it then follows that
(\ref{t3}) is satisfied when both sides act on $F(u)$.
The other commutators in (\ref{45}) (with $\be=0$, $\ga=1$) can now be proved
by using (\ref{47}), (\ref{4b0}) if $B_0$ is involved, and by letting
both sides act on $F(u)$ if $C_0$ is involved.
\eof

\bigbreak
We want to find a realisation of our QISM I algebra as in Theorem \ref{t6}.
{}From this point of view the transformations of $A_0,B_0$, etc.\ as given in
the last statement of the Lemma \ref{t9} do not mean any essential change.
Thus, without lack of generality we may assume that
$\ga=1$ and we may restrict our attention to three special choices
for the pair $\al,\de$: one with $\al,\de\ne0$, one with
$\al\ne0=\de$ and one with $\al=\de=0$.

We now make the restrictive assumption that
$A_0$ and $D_0$ are first order differential operators of the
form (\ref{t56}), analytic on a certain region:
\begin{equation}
A_0=A_{00}(x)+A_{01}(x)\partial_x,\quad D_0=D_{00}(x)+
D_{01}(x)\partial_x.
\label{4xx}\end{equation}
We want to classify solutions of equations (\ref{47}), (\ref{4b0}), (\ref{t3})
such that $A_0$ and $D_0$ have the form (\ref{4xx}).
Note that (\ref{47}) coincides with (\ref{t55}).
Essentially, up to equivalence under gauge transformations and transformations
of the $x$-variable,
all operators $A_0$ and $D_0$ satisfying (\ref{t55}) are given
in the list (\ref{t73})--(\ref{t77}).
The corresponding $B_0$, which equals $[D_0,A_0]$ by (\ref{47}) and
(\ref{4b0}), is also given there. From (\ref{t3}) there 
follows now an expression for $B_0C_0$.
It depends yet on the constant $Q_0$.
If $\al\de$ and $Q_1$ are not both zero (all but the last case)
then we may fix $Q_0$, after a
possible translation of $u$, such that
$\De(u_0)=0$ for a certain $u_0$.

There are now two methods to proceed.
The first method tries to obtain $C_0$ from the known expression for
$B_0C_0$, which may involve
the taking of the inverse of a first order differential operator.
The second method uses Proposition \ref{t90}.
For an explicit pair of shift operator actions on functions $F(u)$
as in Section 3,
we may verify the assumptions of that Proposition. 
Next we define $C_0$ by $C_0F(u)=-uF(u)$.
Sometimes it is not evident that the functions $F(u)$
are linearly independent. Without this property, it is of course not possible
to define $C_0$ as in Proposition \ref{t90}.

Below we follow the first method first in a formal way.
We will give, parallel to the list (\ref{t73})--(\ref{t77}), formal
expressions for $C_0$, or rather for $C(u)=C_0+u$,
such that (\ref{t3}) is satisfied.
Afterwards we will specify a space $W$ on which the operators act, such
that the formal inverse can be understood rigorously.
The expressions for $C(u)$ are as follows.
\begin{eqnarray}
&{\rm Type}\;(A)&\quad
C(u)= -(x\pa_x+b)^{-1}\,(x(1-x)\pa_x^2\nonumber\\
&&\qquad\qquad +(c-(b+a+u+1)x)\pa_x-ab-bu).\label{t79}\\
&{\rm Type}\;(B)&\quad
C(u)=-(x\pa_x^2+(c-x)\pa_x-a-u).\label{t80}\\
&{\rm Type}\;(C')&\quad
C(u)=(\pa_x-1)^{-1}\,(x\pa_x^2+(c+u-x)\pa_x-u).\label{t81}\\
&{\rm Type}\;(D')&\quad
C(u)=\pa_x^2-x\pa_x+u.\label{t82}\\
&{\rm Type}\;(C'')&\quad
C(u)=\tfrac12\pa_x^{-1}\,(x\pa_x^2+(2u+1)\pa_x-x),\quad Q_0=1.\label{t83}
\end{eqnarray}

For each of the five non-trivial types above we will now give a space $W$
on which the operators $A_0,B_0,C_0,D_0$ act, such that $B_0$ is injective,
the inverses of differential operators for types $A$, $C'$ and $C''$
can be rigorously understood, and 
the conditions of Theorem \ref{t6} are satisfied. For this last task 
we have to
give suitable subsets $\FSD$ of $\CC$ such that the equation
$C(u) w=0$ has one-dimensional solution in $W$ for each $u\in\FSD$.

\medbreak\noindent
{\it Type $A$.}\quad
See (\ref{t73}), (\ref{t79}), (\ref{t66}), (\ref{314}),
(\ref{315})--(\ref{316}).
Assume that $b,c\ne0,-1,-2,\ldots\;$.
Let $W$ be the set of all analytic functions on the open unit disk in $\CC$.
Then $A_0, B_0, D_0$ act on $W$, the operator
$B_0$ is moreover injective on $W$ and
$(x\pa_x+b)^{-1}$ acts on convergent power series by termwise application
according to the rule
$$
(x\pa_x+b)^{-1}\,(x^k)=(b+k)^{-1}\,x^k,\quad k=0,1,2,\ldots\;.
$$
so $C(u)$ acts on $W$.
Thus the three equivalent conditions of Lemma \ref{t9} are satisfied.
Then the conditions of Theorem \ref{t6} are satisfied with
$F(u)(x)={}_2F_1(a+u,b;c;x)$
and we can take $\FSD=u_0+\ZZ$ if $u_0\notin(-a+\ZZ)\cup(c-a+\ZZ)$
or we can take $\FSD=\{-a,-a-1,\ldots\}$, for which $F(u)$ becomes a
Jacobi polynomial.

\medbreak\noindent
{\it Type $B$.}\quad
See (\ref{t74}), (\ref{t80}), (\ref{t69}), (\ref{31}),
(\ref{32})--(\ref{33}).
Assume that $c\ne0,-1,\ldots\;$.
Let $W$ be the set of all entire analytic functions on $\CC$.
Then the conditions of Theorem \ref{t6} are satisfied with
$F(u)(x)={}_1F_1(a+u;c;x)$
and we can take $\FSD=u_0+\ZZ$ if $u_0\notin(-a+\ZZ)\cup(c-a+\ZZ)$,
or we can take $\FSD=\{-a,-a-1,-a-2,\ldots\}$, for which $F(u)$ becomes a
Laguerre polynomial.

\medbreak\noindent
{\it Type $C'$.}\quad
See (\ref{t75}), (\ref{t81}), (\ref{31}),
(\ref{351})--(\ref{352}).
Let $W$ consist of all analytic functions $f$ on $(0,\iy)$ such that
$f(x)$ and all its derivatives $f^{(p)}(x)$ are of at most polynomial
growth as $x\to\iy$.
For $f\in W$ define
$$
(\pa_x-1)^{-1}\,f(x)=-e^x\int_x^\iy e^{-y}\,f(y)\,dy.
$$
Then $C(u)$ acts on $W$.
Then the conditions of Theorem \ref{t6} are satisfied with
$F(u)(x)=\Psi(u,c+u;x)$
and we can take $\FSD=u_0+\ZZ$ if $u_0\notin\ZZ$,
or we can take $\FSD=\{0,-1,-2,\ldots\}$, for which $F(u)$ becomes a
Laguerre polynomial.

\medbreak\noindent
{\it Type $D'$.}\quad
See (\ref{t76}), (\ref{t82}), (\ref{t70}),(\ref{t71}),
(\ref{t63})--(\ref{t62}).
Let $W$ be the set of all entire analytic functions $f$ on $\CC$
such that $f(x)$ and all its derivatives $f^{(p)}(x)$ are on $(0,\iy)$
of at most polynomial growth as $x\to\iy$.
Then the conditions of Theorem \ref{t6} are satisfied with
$F(u)(x)=e^{\frac14 x^2} D_{-u}(x)$
and we can take $\FSD=u_0+\ZZ$ if $u_0\notin\ZZ$,
or we can take $\FSD=\{0,-1,-2,\ldots\}$, for which $F(u)$ becomes a
Hermite polynomial.

\medbreak\noindent
{\it Type $C''$.}\quad
See (\ref{t77}), (\ref{t83}), (\ref{355}),
(\ref{t65})--(\ref{356}).
Let $W$ consist of all analytic functions $f$ on $(0,\iy)$ such that
$f(x)$ and all its derivatives $f^{(p)}(x)$ tend to 0 faster
than any inverse power of $x$ as $x\to\iy$.
For $f\in W$ define
$$
\pa_x^{-1}\,f(x)=-\int_x^\iy f(y)\,dy.
$$
Then $C(u)$ acts on $W$.
Then the conditions of Theorem \ref{t6} are satisfied with
$F(u)(x)=x^{-u} K_u(x)$
and we can take $\FSD=u_0+\ZZ$ for any $u_0$.

\begin{remark}
In a sense, the Type $A$ case is the generic case, since the other cases
can be obtained from it by suitable limit transitions. This can be seen on
the level of formulas for the special functions (cf.\ (\ref{t99})),
of Lie algebras and of QISM I algebra representations.

In \S3.1 we gave several other shift operator pairs of Type $A$.
The pair (\ref{vv1})--(\ref{vv2}) is a variant of
the pair (\ref{315})--(\ref{316}). Indeed, first replace
$F(a+u,b;c;x)$ by a second solution to
(\ref{314}) as in (\ref{t100}), then make a transformation as in
(\ref{t101}), next a transformation of the independent variable, and finally a
gauge transformation.
If $a=0$ and $c\notin\ZZ$ in (\ref{vv1})--(\ref{vv2}) then we can take for
$W$ the set of all polynomials in $(1-x)^{-1}$ and for $\FSD$ the set
$\{0,-1,-2,\ldots\}$.
For $u\in\FSD$ the $F(u)(x)=(1-x)^u\,F(u,b+u;c+u;x)$ is
a polynomial of exact degree $-u$
in $(1-x)^{-1}$. Thus the functions $F(u)$ are linearly independent elements
of $W$ and we realize on the span of the $F(u)$ a Type $A$ representation
of the QISM I algebra because of Proposition \ref{t90}.

The shift operator pair (\ref{t95})--(\ref{t96}) can be obtained by
specialization of (\ref{vv1})--(\ref{vv2}) (for negative $k$ also apply a
gauge transformation). Superficially one would say that equations
(\ref{t95})--(\ref{t96}) realize a finite-dimensional representation
of the QISM I algebra on the span of the $F_k$ for $k=-n,\ldots,n$.
However, $(1-x)^n\,F_k(x)$ is a polynomial in $x$ of exact degree $n$.
So the $2n+1$ functions $F_k$ can never be linearly independent and it is
impossible to have an operator $C_0$ with $C_0\,F_k=-kF_k$.
\end{remark}

\section{Rank 1 $L$-operators for the QISM II algebra}
\setcounter{equation}{0}
The (unitary) $U$-algebra is the algebra with the matrix elements of $U(u)$
as generators and with relations (\ref{4}) and (\ref{141}).
Let us pass to a quotient algebra by adding relations stating that
$(u-\frac12)U(u)$ is a polynomial of degree $\le3$ in $u$.
In other words, we make the ansatz that $U(u)$ is of the form
\begin{eqnarray}
{L(u)}&=&\left(\matrix{A(u)&B(u)\cr C(u)&D(u)}
\right)\nonumber\\
&=&(u-\tfrac12)^2L_2+
(u-\tfrac12)L_1+L_0+
(u-\tfrac12)^{-1}L_{-1},\label{51}\\
\noalign{\hbox{where}}
{L_2}&=&\left(\matrix{A_2&B_2\cr C_2&-A_2}
\right),\quad L_1=\left(\matrix{A_1& B_2\cr 
 C_2&A_1-2 A_2}\right),\nonumber\\
{L_0}&=&\left(\matrix{A_0&B_0+\tfrac14B_2\cr 
C_0+\tfrac14C_2&-A_0+2 A_1-2A_2}
\right),\quad L_{-1}=\left(\matrix{A_{-1}&0\cr 
0&A_{-1}}\right).\label{51-xx}\end{eqnarray}
Substitute (\ref{51}) in relation (\ref{4}).
Then we get the algebra with the matrix elements of the $L_i$ ($i=2,1,0,-1$)
as generators and with relations
\begin{eqnarray}
{L_2^{(1)}L_{2,1,0,-1}^{(2)}}&=&L_{2,1,0,-1}^{(2)}L_2^{(1)},\quad 
L_{-1}^{(1)}L_{2,1,0,-1}^{(2)}=L_{2,1,0,-1}^{(2)}L_{-1}^{(1)},
\label{52}\\
{[L_1^{(1)},L_0^{(2)}]}&=&- [P,L_2^{(1)}L_0^{(2)}]-
 L_2^{(1)}P L_0^{(2)}+ L_0^{(2)}P L_2^{(1)},
\label{53}\\
{[L_0^{(1)},L_0^{(2)}]}&=&- [\{P,L_1^{(1)}\},L_0^{(2)}]-
2 A_{-1}[P,L_2^{(1)}]+{[L_0,L_2]}^{(2)}.
\label{54}\end{eqnarray}
Here curved brackets mean anticommutator.
The relations (\ref{52}) imply that the entries of the $L_2$ and $L_{-1}$
matrices 
are in the center of the algebra. Let us pass once more to a quotient algebra
by adding the relations
\begin{equation}
L_2=\left(\matrix{\alpha&\beta\cr \gamma&-\alpha}\right),
\quad A_{-1}=\delta,
\label{55}\end{equation}
for certain $\al,\be,\ga,\de\in\CC$.
We thus obtain an algebra with generators $A_0,A_1,B_0,C_0$ and relations
\begin{eqnarray}
{[A_1,A_0]}&=&\gamma B_0-\beta C_0,\nonumber\\
{[A_1,B_0]}&=&-2\alpha B_0+\beta\left(2A_0-2A_1+\tfrac32\alpha
\right),\nonumber\\
{[A_1,C_0]}&=&2\alpha C_0+\gamma\left(-2A_0+2A_1-\tfrac32\alpha
\right),\nonumber\\
{[A_0,B_0]}&=&-\{A_1,B_0\}+\beta\left(2A_0-\tfrac52 A_1+2\alpha
+2\delta\right),\nonumber\\
{[A_0,C_0]}&=&\{A_1,C_0\}+\gamma\left(-2A_0+\tfrac52 A_1-2\alpha
-2\delta\right),\nonumber\\
{[B_0,C_0]}&=&-2\{A_0,A_1\}+4(A_1-\alpha)^2+4\alpha A_0+4\alpha\delta\,.
\label{56}\end{eqnarray}
The quantum determinant (\ref{10}) now has the form
\begin{eqnarray}
\Delta(u)&=&-({\alpha}^2+\beta\gamma)u^4+Q_2u^2+Q_0+
{\delta}^2\, u^{-2}\,,\label{57}\\
Q_2&=&A_1^2-2\alpha A_0-\gamma B_0-\beta C_0+\tfrac12\beta\gamma,
\label{t20}\\
Q_0&=&-A_0^2-B_0C_0+2\delta A_1-\tfrac14\gamma B_0-\tfrac14\beta C_0-
\tfrac1{16}\beta\gamma.
\label{t21}\end{eqnarray}
Here the right hand sides of (\ref{t20}) and (\ref{t21}) give operators
in the center of the algebra. We consider (\ref{t20}) and (\ref{t21}) as
added relations, for a certain choice of $Q_2,Q_0\in\CC$.
So a representation of the algebra with relations (\ref{52})--(\ref{54})
which has the
property that all elements in the center of the algebra are represented as
scalars, can also be viewed as a representation of the algebra with
relations (\ref{56}), (\ref{t20}) and (\ref{t21})
for a certain choice of $\al,\be,\ga,\de,Q_0,Q_2$.

Let us assume that
\begin{equation}
\alpha=\beta=0,\quad\gamma=1.
\label{t22}
\end{equation}
First we make an observation somewhat analogous to Remark \ref{t102}.

\begin{remark}\label{t103}
Consider the algebra $\FSA$ with generators $A_0,A_1,B_0,C_0,\de$ and with
two sets of relations: relations (\ref{56}) under assumption (\ref{t22}), and
relations stating that $\de$ is in the center of the algebra.
There is a homomorphism of this algebra into the universal enveloping
algebra $\FSU({\rm e}(3))$
of the Lie algebra e(3). A Lie group corresponding to $e(3)$ is the
group of motions of 3-dimensional Euclidean space.
The Lie algebra e(3) is 6-dimensional.
It can be described by a basis
${\cal P}^\pm,{\cal P}^3,{\cal J}^\pm,{\cal J}^3$ and
commutation relations
$$
[{\cal J}^3,{\cal J}^\pm]=\pm{\cal J}^\pm,\qquad
[{\cal J}^3,{\cal P}^\pm]=[{\cal P}^3,{\cal J}^\pm]=\pm{\cal P}^\pm,
$$
$$
[{\cal J}^+,{\cal P}^+]=[{\cal J}^-,{\cal P}^-]=[{\cal J}^3,{\cal P}^3]=0,
$$ 
$$
[{\cal J}^+,{\cal J}^-]=2{\cal J}^3,\qquad
[{\cal J}^+,{\cal P}^-]=[{\cal P}^+,{\cal J}^-]=2{\cal P}^3,
$$
$$
[{\cal P}^3,{\cal P}^\pm]=[{\cal P}^+,{\cal P}^-]=0.
$$ 
The center of the universal enveloping algebra is generated by two
Casimir elements:
$$
C={\left({\cal P}^3\right)}^2+{\cal P}^+{\cal P}^-,
\qquad {\tilde C}=\tfrac12\left({\cal P}^+{\cal J}^-+
{\cal P}^-{\cal J}^+\right)+{\cal P}^3{\cal J}^3.
$$

It is now straightforward to verify that the relations for the generators of
$\FSA$ are satisfied when we put these generators equal to the following
elements of $\FSU({\rm e}(3))$.
\begin{eqnarray}
A_0=\tfrac12({\cal P}^+{\cal J}^--{\cal P}^-{\cal J}^+),\quad
A_1={\cal P}^3,\quad B_0=-{\cal P}^+{\cal P}^-,
\nonumber\\
C_0=-\tfrac12\{{\cal J}^+,{\cal J}^-\}-{\left({\cal J}^3\right)}^2
-\tfrac14,\quad
\de=-{\tilde C}{\cal J}^3.\qquad\nonumber
\end{eqnarray}
This yields the announced algebra homomorphism of $\FSA$ into
$\FSU({\rm e}(3))$. We do not yet know if this homomorphism is injective.
\end{remark}

The following lemma can be proved in a straightforward way. It
shows that equations (\ref{56}), (\ref{t20}) and (\ref{t21}), with
(\ref{t22}), and under the assumption that $B_0$ is injective,
can be equivalently written in a much more simple form.

\begin{lemma}\label{t23}
Let $\delta,Q_0,Q_2$ be scalars.
Let $A_0,A_1,B_0,C_0$
be operators acting on some linear space. Let $B_0$ be injective.
Then the following three statements are equivalent:
\vskip 2mm

(a)\quad
$\pmatrix{(u-\tfrac12)A_1+A_0+\de(u-\tfrac12)^{-1}&B_0\cr
u^2+C_0&(u+\tfrac32)A_1-A_0+\de(u-\tfrac12)^{-1}}$
\vskip 2mm

is a representation of the QISM II algebra with quantum determinant 

$\De(u)=Q_2u^2+Q_0+\de^2u^{-2}$;

(b)\quad
The six commutators (\ref{56}) and formulas (\ref{t20}), (\ref{t21})
are valid with $\al=\be=0$

and $\ga=1$;

(c)\quad
The following three equations are valid:
\begin{eqnarray}
[A_1,A_0]&=&A_1^2-Q_2,\label{59}\\
{B_0}&=&A_1^2-Q_2,\label{58}\\
B_0C_0&=&2\de A_1-A_0^2-\tfrac14  B_0-Q_0.\label{t24}
\end{eqnarray}
Moreover, if
$\{A_0,A_1,B_0,C_0,\de,Q_0,Q_2\}$ satisfy these equivalent conditions
then so do
$\{\la A_0,\la A_1,\la^2 B_0,C_0,\la\de,\la^2 Q_0,
\la^2 Q_2\}$
(where $\la$ is a nonzero scalar).
\end{lemma}

\bigbreak
We want to find a realisation of our QISM II algebra as in Theorem \ref{t7}.
{}From this point of view the transformations of $A_0,A_1$, etc.\ as given in
the last statement of the Lemma \ref{t23} do not mean any essential change.
Thus, without lack of generality
we may restrict our attention to two special choices
for the $Q_2$: one with $Q_2\neq 0$ and one with
$Q_2=0$.

We now make the restrictive assumption
that $A_0$ is a first order differential
operator and $A_1$ is a scalar function of $x$:
\begin{equation}
A_0=A_{00}(x)+A_{01}(x)\partial_x,\quad A_1=A_{10}(x).
\label{510}
\end{equation}
The following approach should now be followed.
Find all operators of the form (\ref{510}) such that (\ref{59}) is satisfied
for some number $Q_2$.
(It is sufficient to find one solution in each equivalence class formed
by gauge transformations and transformations of the $x$-variable.)$\;$
Then define $B_0$ by (\ref{58}) and try to define $C_0$ by (\ref{t24}).
Fix some function space $W$ on which these operators act.
Then the equivalent conditions of Lemma \ref{t23} are satisfied.
Finally check if the conditions of Theorem \ref{t7} are satisfied for some
choice of $\FSD$.
Note that the analogue of the second method described in \S5 cannot be used
here, since we were not able to formulate an analogue for the QISM II case of
Proposition \ref{t90}.

For the case $\de=0$ a classification (up to equivalence) of all operators
of the form (\ref{510}) such that (\ref{59}) is satisfied was already given
(in the non-trivial cases) by (\ref{v73}) and (\ref{v75}). Only
types $A$ and $C''$ showed up.
We generalize these results for the case of general $\de$ in the short list
below. It is immediately verified that equation (\ref{59}) is satisfied
for $A_0$ and $A_1$ given there. It turns out that $B_0$ is a function,
so (\ref{t24}) defines $C_0$ without problems.
We can show that the possibilities for $A(u), B_0$ and $C(u)$ listed
below are the only ones up to gauge transformations and transformations
of the $x$-variable, but we do not include the proof here.

\begin{eqnarray}
&&\hbox{Generalized Type ($A$)}
\quad A(u)=(u-\tfrac12)(x-\tfrac12)+x(1-x)\pa_x+
(\tfrac12-a)x\qquad\qquad\nonumber\\
&&\qquad\qquad\qquad\qquad\qquad\qquad
+\tfrac12 c-\tfrac12+{\de\over u-\tfrac12},
\quad B_0=[A_1,A_0]=-x(1-x),\nonumber\\
&&\qquad\qquad\qquad\qquad C(u)=x(1-x)\pa_x^2+[c-(2a+1)x]\pa_x-a^2+u^2;
\label{last1}\\
&&\hbox{Generalized Type ($C''$)}\quad A(u)=(u-\tfrac12)x^{-1}+\pa_x+
\tfrac12 x^{-1}
+{\de\over u-\tfrac12}\,,\qquad\qquad\nonumber\\
&&\qquad\qquad B_0=[A_1,A_0]=x^{-2},
\quad C(u)=-x^2\pa_x^2-x\pa_x-x^2+2\de x+u^2.
\label{last2}\end{eqnarray}
For these two $L$-operators we can give functions $F(u)$ on which 
$A(u)$ and $-A(-u)$ act as shifting operators. See equations
(\ref{30001})--(\ref{30002}), (\ref{last3})--(\ref{last4}),
respectively. Below we give a space $W$ and
suitable subsets $\FSD$ of $\CC$ such that the equation
$C(u) w=0$ has one-dimensional solution in $W$ for each $u\in\FSD$,
so such that the conditions of Theorem \ref{t7} are satisfied.

\medbreak\noindent
{\it Generalized Type $A$.}\quad
See (\ref{last1}), (\ref{t66}), (\ref{314}),
(\ref{30001})--(\ref{30002}).
Assume that $c\ne0,-1,$ $-2,\ldots\;$.
Let $W$ be the set of all analytic functions on the open unit disk in $\CC$.
Then $A(u)$, $-A(-u)$, $B_0$, $C(u)$ act on $W$.
Then the conditions of Theorem \ref{t7} are satisfied with
$F(u)(x)={}_2F_1(a+u,a-u;c;x)$
and we can take $\FSD=u_0+\ZZ$ if $u_0\notin(\pm a+\ZZ)\cup(\pm(c-a)+\ZZ)$,
or we can take $\FSD=\{-a,-a-1,\ldots\}$
if $-2a,-2a+c\ne 0,1,2,\ldots\,$,
or we can take $\FSD=\{a,a+1,\ldots,-a\}$
if $a=0,-1,-2,\ldots\,$.
In the second and third case $F(u)$ becomes a Jacobi polynomial with fixed
parameters, while the shift only affects the degree.

\medbreak\noindent
{\it Generalized Type $C''$.}\quad
See (\ref{last2}), (\ref{31}), (\ref{last3})--(\ref{last4}).
Let $W$ consist of all analytic functions $f$ on $(0,\iy)$ such that
$f(x)$ and all its derivatives $f^{(p)}(x)$ tend to 0 faster
than any inverse power of $x$ as $x\to\iy$.
Then $A(u),-A(-u),B_0, C(u)$ (see (\ref{last2})) act on $W$.
Then the conditions of Theorem \ref{t7} are satisfied with
$F(u)(x)=x^ue^{-x/2}\Psi(2\de+\tfrac12+u,2u+1;x)$
and we can take $\FSD=u_0+\ZZ$ if $u_0\notin(\pm(2\de+{1\over 2})+\ZZ)$,
or we can take $\FSD=\{-2\de-\tfrac12,-2\de-\tfrac12-1,\ldots\}$
if $-4\de-1\ne 0,1,2,\ldots\,$,
or we can take $\FSD=\{2\de+\tfrac12,2\de+\tfrac12+1,\ldots,
-2\de-\tfrac12\}$ if $2\de+\tfrac12=0,-1,-2,\ldots\,$.
In the second and third case $F(u)$ becomes a 
Laguerre polynomial in $x^{-1}$ multiplied by an exponential and a power,
while the shift affects both the degree and the parameter.

\begin{remark}
Write the operator in the left hand side of (\ref{v2}) as
$A(u)=(u-\tfrac12)A_1+A_0$. Then equation (\ref{59}) is satisfied.
In fact, this operator is equivalent to the Type $(A)$ case given in
(\ref{v73}). However, the Legendre functions on which the shift operator
pair in (\ref{v2})--(\ref{v3}) acts, cannot be obtained generally from the
hypergeometric functions in the $\de=0$ case of
(\ref{30001})--(\ref{30002}) by just making a gauge transformation and
a change of $x$-variable. We have to pass also to another solution of the
corresponding second order differential equation.
The choice of an appropriate space $W$ is not so clear now.
But in the case of a finite dimensional representation we can pass
{}from (\ref{30001})--(\ref{30002}) to (\ref{v2})--(\ref{v3}) without passing
to another solution of the differential equation.

For the case $\de=0$ of (\ref{last3})--(\ref{last4}) the operators
coincide with the operators in (\ref{354})--(\ref{t64}).
However, the functions in (\ref{last3})--(\ref{last4}) do not specialize
for $\de=0$ to the Bessel functions in (\ref{354})--(\ref{t64})
but to other solutions of the corresponding second order differential
equation. For the functions in (\ref{354})--(\ref{t64}) there may be a problem
of a good choice of $W$.

Note that, in the cases of a finite dimensional representation of the
QISM II algebra we met above, the functions $F(u)$ and $F(-u)$
($u\in\FSD$) are proportional and certainly not linearly independent.
This is compatible with the fact that $F(u)$ must be eigenfunction of $C_0$
with eigenvalue $-u^2$.
\end{remark}

\begin{remark} 
Let us give more comments on the homomorphism
of the quadratic algebra $\FSA$ coming from the QISM II
algebra into $\FSU$(e(3)) (see Remark \ref{t103}). 
Miller defines \cite{mi68} the following two operators 
in the universal enveloping algebra of the algebra e(3):
$$
X(u,+)={\cal P}^-{\cal J}^++{\cal P}^3{\cal J}^3+
(u+1){\cal P}^3-{{\tilde C}\over u+1}{\cal J}^3-{\tilde C},
$$
$$
X(u,-)=-{\cal P}^-{\cal J}^+-{\cal P}^3{\cal J}^3+
u{\cal P}^3-{{\tilde C}\over u}{\cal J}^3+{\tilde C},
$$
Then he gets the following actions for these operators on the 
basis vectors $f_m^{(u)}$ of the representation space for the
algebra e(3):
$$
X(u,+)f_m^{(u)}={\omega(u-q+1)\over u+1}f_m^{(u+1)},
$$
$$
X(u,-)f_m^{(u)}={\omega(u+m)(u-m)(u+q)\over u}f_m^{(u-1)},
$$
where the basis functions $f_m^{(u)}$ are fixed by the diagonal action of the 
following four mutually commuting operators:
$$
{\cal J}^3f_m^{(u)}=mf_m^{(u)},\qquad \left(
\tfrac12\{{\cal J}^+,{\cal J}^-\}+{\left({\cal J}^3\right)}^2
\right)f_m^{(u)}=u(u+1)f_m^{(u)},
$$
$$
Cf_m^{(u)}=\omega^2f_m^{(u)},\qquad {\tilde C}f_m^{(u)}=\omega qf_m^{(u)}.
$$
Notice that the operators $X(u,\pm)$ shift the parameter $u$ 
of the basis (while the ${\cal J}^\pm$ shift $m$ and constitute, together 
with ${\cal J}^3$, the sub-algebra sl(2)).
The operators $X(u-\tfrac12,\pm)$ 
coincide with our operators $\mp A(\mp u)=(u\pm\tfrac12)A_1\mp A_0+
\de(u\pm\tfrac12)^{-1}$ and they give some shift operator
actions for Gauss and confluent hypergeometric functions
(see generalized Types $A$ and $C''$ (6.16)--(6.17), which are 
Types $E$ and $F$, respectively, in \cite{mi68}).
It would be natural to try to find any analogous
homomorphism of the algebra $\FSA$ with  general commutation 
relations (6.7) (when $\al$ and $\be\neq 0$)
into the universal enveloping algebra of a Lie algebra.
We found that it is possible to do so with $\FSU$(o(4))
in the  case $\be\neq 0$,
$\al=0$.
Then the case $\be=0$ corresponds to the
contraction of o(4) to e(3). The case of $\al\neq 0$ is still unsolved.
\end{remark}

\section{Concluding remarks}
\setcounter{equation}{0}
In this Section we would like to give some comments
on the possible applications of the above results in the 
theory of finite-dimensional quantum integrable systems.
By use of the comultiplication operation we get the monodromy
matrix for the quantum integrable chain as a product of
$L$-operators each being associated with a particular 
site of a chain.
The main question is to study the spectral problem for 
the complete set of commuting integrals of motion.
In such a way the special functions appear as common 
eigenfunctions of those commuting operators. 

In the present paper we have constructed a lot of new
$L$-operators for both $T$- and $U$-algebras connecting
each particular $L$-operator with a particular recurrence
relation for the corresponding special function. In this approach 
we have got an interpretation for the spectral parameter $u$
appearing as an argument of operator-valued entries of
the matrix $L(u)$. The meaning of the spectral parameter $u$
is that it is a parameter (like $a$, $b$, or $c$ in 
${}_2F_1(a,b,c;x)$) of the special function $F(u)$ defined by 
the following rule:
\begin{equation}
C(u)\,F(u)=0.
\label{61}\end{equation}
This equation looks like one appearing in the so-called
algebraic Bethe ansatz (ABA) technique \cite{ks82,sk91,tf79}. 
In this
analogy the $F(u)$ is a pseudovacuum state. But the
crucial difference now is that our ``pseudovacuum'' might 
depend on the spectral parameter, so it can exist in
the situations where the standard ABA does not work. 
The method of variable separation or functional Bethe
ansatz (FBA)
\cite{ku89,ku90,sk84,sk91} was developed to overcome
the obstacles in application of ABA. This method
has dealt with operator zeros of the equation
\begin{equation}
C(u)=0.
\nonumber\end{equation}
For $F(u)$ defined by (\ref{61})
we have been able to find the operators acting
in the spectral parameter space (see Theorems \ref{t6} and \ref{t7}).
These operators act as shifting operators to the function
$F(u)$ considered now as a function of spectral parameter.
We have work in progress on
further applications of this idea.
In particular, we are preparing a paper with a
further generalisation of this technique
for $q$-special functions.
\pagebreak

\end{document}